\documentclass{elsart}

\journal{Linear Algebra and its Applications}

\usepackage{ifpdf}
\ifpdf
\usepackage[pdftex]{graphicx}
\else
\usepackage[dvips]{graphicx}
\fi

\usepackage{amsmath,amssymb}
\usepackage{mathdots}
\usepackage{url}

\newcommand{\p}[1]{[#1]}
\newcommand{\plucker}[4]{\frac{\p{#1}\p{#2}}{\p{#3}\p{#4}}}
\newcommand{\iplucker}[4]{\p{#1}\p{#2} / \p{#3}\p{#4}}

\newcommand{\R}{\mathbb{R}}

\begin{document}

\begin{frontmatter}

\title{On generators of bounded ratios of minors for totally positive matrices}

\date{13 October 2007}

\author[ND]{Adam Boocher\corauthref{cor}},
\corauth[cor]{Corresponding author.}
\ead{aboocher@nd.edu}
\author[UMN]{Bradley Froehle}
\ead{bfroehle@math.berkeley.edu}

\address[ND]{University of Notre Dame, Department of Mathematics\\
Notre Dame, Indiana 46556}
\address[UMN]{University of California, Department of Mathematics\\
Berkeley, CA 94720, United States}

\begin{abstract}
We provide a method for factoring all bounded ratios of the form
$$\det A(I_1|I_1')\det A(I_2|I_2')/\det A(J_1|J_1')\det A(J_2|J_2')$$
where $A$ is a totally positive matrix, into a product of more
elementary ratios each of which is bounded by 1, thus giving a new
proof of Skandera's result.  The approach we use generalizes the
one employed by Fallat et al.~in their work on principal minors.
We also obtain a new necessary condition for a ratio to be bounded
for the case of non-principal minors.
\end{abstract}


\end{frontmatter}

\section{Introduction}

An $n\times n$ matrix $A$ is called totally positive if every
minor of $A$ is positive.
If $I,I' \subseteq \{1,2,\ldots,n\}$ with
$|I|=|I'|$, we denote the minor of $A$ with row set $I$ and column
set $I'$ as $(I|I')(A):=\det A(I|I')$.  If
$S=((I_1|I_1'),\ldots,(I_p|I_p'))$ is a sequence of $p$ row and
column sets, we define a function $S(A)=\det A(I_1|I_1')\cdot\det
A(I_2|I_2')\cdots \det A(I_p|I_p')$. 
Please note that
$S(A)>0$ for any choice of $S$ and for all totally positive matrices $A$.
Similarly, if $T=((J_1|J_1'),\ldots,(J_q|J_q'))$ is another sequence
of $q$ row and column sets, we say that $S\leq T$ (with respect to
the class of totally positive matrices) if $S(A)\leq T(A)$ for all totally positive matrices
$A$. Note that if we take the convention that
$(\emptyset|\emptyset)(A)=1$, we are free to assume that $S$ and $T$
are both sequences of the same size (i.e., $p=q$) by appending an
appropriate number of $(\emptyset|\emptyset)$ to the shorter
sequence.

It is also reasonable to ask when the ratio $S(A)/T(A)$ is bounded
by some $k>0$ for all totally positive matrices $A$.  If this is
true, we say that the ratio $S/T$ is \emph{bounded by $k$}.  It is
clear that $S\leq T$ if and only if $S/T$ is bounded by $1$. It has
been conjectured that if $S/T$ is bounded (by any number), then it
is necessarily bounded by $1$ (e.g., see \cite{Fallat:2000ps}).

Recently, the problem of classifying all such ratios and
inequalities has been a subject of much interest.  Fallat et al.~\cite{Fallat:2003az} were able to classify a large class of ratios
of products of principal minors.  In particular, they gave necessary
and sufficient conditions for a ratio of products of two minors to
be bounded over totally positive matrices.  This result was later generalized to
the case of non-principal minors by Skandera \cite{Skandera:2004qp}.
In this paper we generalize a necessary condition in
\cite{Fallat:2003az} to the case of non-principal minors, and our
main result is an explicit factorization of ratios of the form
\[\frac{\det A(I_1|I_1') \det A(I_2|I_2')}{\det A(J_1|J_1') \det A(J_2|J_2')}\]
into products of elementary ratios.  This in particular implies the
result of Skandera describing bounded ratios of this form.  It has
been conjectured by Gekhtman that all bounded ratios are products
of these elementary ratios \cite{Fallat:2003az}.

 \subsection{Planar Networks and Totally Positive Matrices}\label{subsec:planar_networks}
The relationship between totally positive matrices and directed
acyclic weighted planar networks is well studied.  It was first discussed by Karlin and
McGregor in 1959 \cite{Karlin:1959mc}.  For a more modern presentation, refer to 
the paper by Fomin and Zelevinsky \cite{Fomin:2000ab}.  In an attempt to keep the manuscript mostly self-contained, we will present some relevant results from these papers.

\begin{figure}[htb]
\centering
\includegraphics[width=5.4in]{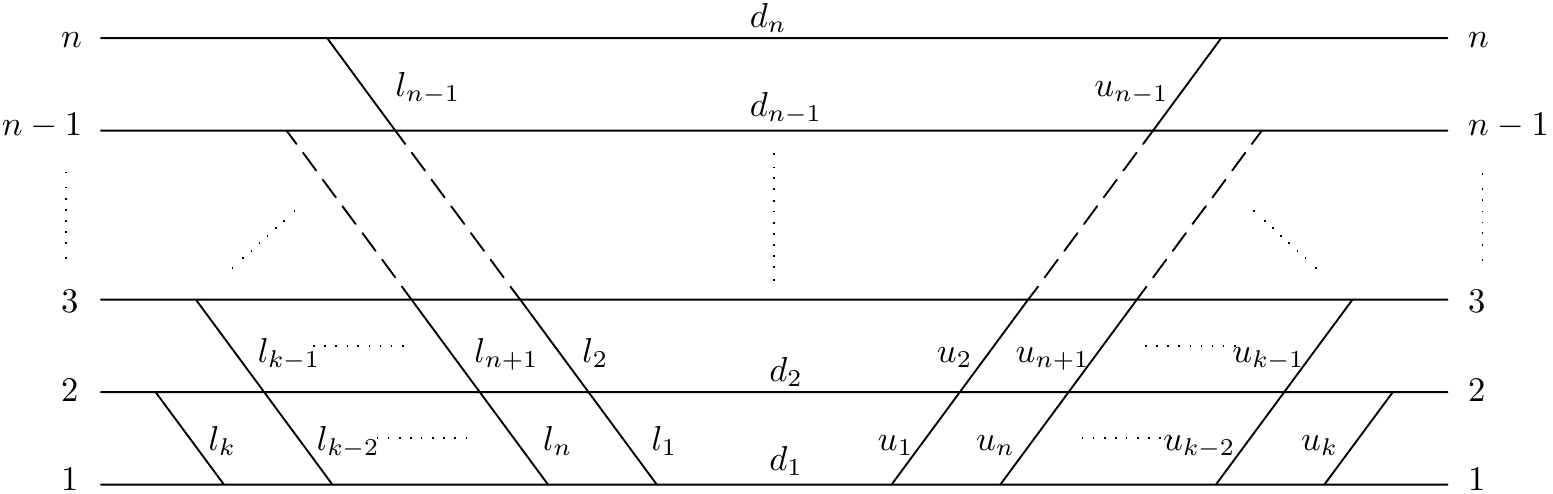}
\caption{\label{fig:general}General Planar Network}
\end{figure}

A typical directed acyclic weighted planar network is shown in
Figure~\ref{fig:general}. Note that because the graph is acyclic, we
can stretch the network in an appropriate fashion so that the
direction of each edge is oriented from left-to-right.  Furthermore,
the network is assumed to have $n$ labeled sources (on the left) and
$n$ labeled sinks (on the right).   Both sources and sinks are
labeled bottom to top.  Additionally, to each edge of the network we
associate a positive weight.  In Figure~\ref{fig:general}, these
weights are shown as $l_i$, $d_j$, or $u_k$.  Unmarked weights are
assumed to be $1$.

Let $\pi$ be any path running left-to-right from source $i$ to sink $j$. We define the weight of this path
to be the product of the weights along each edge of the path and denote this as $w(\pi)$.

To each such diagram, we can associate a totally positive matrix $A$
with entries $a_{ij}$ given by \begin{equation}a_{ij}=\sum_{\pi:i\to
j} w(\pi)\label{bijection} \end{equation} where the summation is
over all paths $\pi$ that begin at source $i$ and end at sink $j$.
Formula (\ref{bijection}) establishes a bijection between totally
positive matrices and planar networks of the kind depicted in
Figure~\ref{fig:general}.  This fact is equivalent to Anne Whitney's
Reduction Theorem \cite{Whitney:1952ab}.

Let us define a path family $\pi$ as a set of non-intersecting paths
running from left-to-right starting at the sources in $I$ and
terminating at the sinks in $I'$.  The weight of such a path family
$w(\pi)$ is defined to be the product of the weights of each path in
the path family.  As shown in \cite{Karlin:1968wh} the minor with
row set $I$ and column set $I'$ is
\[ \det A(I|I') = \sum_{\pi:I\to I'} w(\pi) \]
where the summation is over all such possible path families from $I$
to $I'$.

Given the row set $I$ and column set $I'$, we have found it helpful
to follow Skandera \cite{Skandera:2004qp} in defining the set $I''$
which encapsulates both $I$ and $I'$
\begin{equation}\label{eq:i''}
I'' = I\cup \{2n+1-i \, | \, i \in I'^c \} \end{equation} where
$I'^c = \{1,2,\ldots,n\} \setminus I'$.  While this $I''$ may seem
cryptic, it has a natural interpretation if one considers an
embedding of totally positive matrices into the totally positive
part of the Grassmannian $Gr(n,2n)$.

\subsection{Grassmannians}
\label{sec:grassmannian} In this section we will discuss the
real Grassmannian and refer the reader to Section 5.4 of
\cite{Smith:2000ab} for more information.  Recall that the real Grassmannian $Gr(n,2n)$ is the set of 
$n$-dimensional subspaces of $\R^{2n}$, i.e.
\[
	Gr(n,2n) = \{ \text{Real $2n \times n$ matrices of rank $n$} \} / GL(n,\R)
\]
where we have factored out the action of right multiplication by an invertible $n\times n$ matrix.


It is clear that an element $\Lambda \in Gr(n,2n)$ does not have a unique matrix representation, but rather a collection of matrix representatives which are unique up to right multiplication by an invertible $n\times n$ matrix.

If $A$ is a $2n\times n$ matrix representative of $\Lambda$, we can define the \emph{Pl\"ucker coordinates of $\Lambda$ with respect to $A$} (or more briefly the \emph{Pl\"ucker coordinates of $A$}) to be the vector of all $n\times n$ minors of the matrix $A$, i.e. an element of real ${2n \choose n}$-space.


We say an element $\Lambda \in Gr(n,2n)$ is totally positive if there exists a matrix representative $A$ of $\Lambda$ such that every Pl\"ucker coordinate of $A$ is positive.  The totally positive part of the Grassmannian $Gr(n,2n)$ is then defined to be
\[
	TPGr(n,2n) = \{ \Lambda \in Gr(n,2n) : \text{$\Lambda$ is totally positive} \}.
\]
If $\Lambda \in TPGr(n,2n)$ we say that the \emph{standard matrix representative} of $\Lambda$ is the $2n\times n$ matrix representative $\bar A$ with lower $n\times n$ submatrix equal to
\[ \begin{bmatrix}
& & & 1 \\
& & -1 & \\
& \iddots & & \\
\pm 1 & & &
\end{bmatrix}. \]  Note that such a matrix can always be chosen because the lower $n\times n$ block of any matrix representative of $\Lambda$ is always of full rank.

\begin{prop}
	There is a natural bijection:
	\[
		\{\text{Totally positive $n\times n$ matrices}\} \leftrightarrow TPGr(n,2n).
	\]
\end{prop}

\begin{pf}
	Let $\Lambda \in TPGr(n,2n)$, and let $\bar A$ be its standard matrix representative.  We shall denote the upper $n\times n$ submatrix of $\bar A$ as $A$.  Then the relation
	\[ \det A(I|I') = \det \bar A(I''|(1,2,\ldots,n)) \]
	where $I''$ is defined as in equation~\ref{eq:i''}, and 
	the positivity of all Pl\"ucker coordinates of $\bar A$ imply that $A$ is a totally positive matrix.
	
	This same relation allows us to pass from an $n\times n$ totally positive matrix $A$ to an element $\Lambda \in TPGr(n,2n)$ by choosing $\Lambda$ to be the unique element with standard matrix representative having $A$ as the upper $n\times n$ submatrix.
\end{pf}

With this bijection clearly established we will maintain the convention of using $A$ to represent a totally positive $n\times n$ matrix and $\bar A$ as its corresponding standard matrix representative in $TPGr(n,2n)$.

For additional notational convenience, and to distinguish between minors and
Pl\"ucker coordinates, we will designate
index sets representing Pl\"ucker coordinates using Greek letters
and drop the bar notation where its meaning is unambiguous.  That is, for an index set $\alpha_j
\subset \{1,2,\ldots,2n\}$ of size $n$, we define
\begin{equation}[\alpha_j] (A) := \det \bar A (\alpha_j | (1,2,\ldots,n))
.\label{squarebracketdefinition}\end{equation}

Unless stated otherwise, all index sets $\alpha_j$ in the
remainder of the paper will be assumed to be cardinality $n$ subsets
of $\{1,\ldots,2n\}$.

If we have a sequence of index sets $\alpha =
(\alpha_1,\alpha_2, \ldots,\alpha_p)$, we can define the function
$\alpha(A)$ as a product of Pl\"ucker coordinates $$\alpha(A) = \prod_{i=1}^p
[\alpha_i](A)$$ where $A$ is an $n\times n$ totally positive matrix.

If we similarly let $\beta=(\beta_1,\ldots,\beta_q)$ be another
sequence of index sets, we write $\alpha \leq \beta$ (with respect
to totally positive matrices) if $\alpha(A) \leq \beta(A)$ for all
$n \times n$ totally positive matrices $A$.
We say that $\alpha/\beta$ is bounded by $k$ (with respect
to totally positive matrices) if $\alpha(A)/\beta(A)
\leq k$ for all totally positive matrices $A$.  Note that $\alpha
\leq \beta$ is equivalent to saying $\alpha/\beta$ is bounded by
$1$.

By construction we have
$\p{(n+1,n+2,\ldots,2n)}(A)=1$ for all totally positive matrices
$A$, and thus in general we will assume that $\alpha$ and $\beta$ each
contain the same number of index sets.

Lastly, when we say
$$\frac{[\alpha_1]\cdots[\alpha_p]}{[\beta_1]\cdots[\beta_p]}$$
is bounded (resp.~bounded by $k$), we mean that the ratio $\alpha/\beta$ is bounded (resp.~bounded by $k$) 
where $\alpha = (\alpha_1,\ldots,\alpha_p)$ and 
$\beta = (\beta_1,\ldots,\beta_p)$.


\section{Operations Which Preserve Bounded Ratios}

Before proving the main theorem in the work of Fallat et al.~(see \cite{Fallat:2003az}), they developed several operators that preserved bounded ratios of minors, namely what they called the Complement, Reversal, Shift, Insertion, and Deletion operators.  Of these operators, we will provide generalizations of the shift and reversal operators.  The insertion and deletion operators were not generalized because they have little applicability to our situation in which the cardinality of each index set must remain fixed.  The complement operator was not studied.


\begin{defn}[Cyclic Shift]
\label{defn:cyclicshift} For an index set $\alpha_j$, define a
\emph{cyclic shift} of the elements of $\alpha_j$ as
\[ \sigma(\alpha_j) = \{ i+1 \mod{2n}\, | \, i\in \alpha_j\} \] which
maps $i\in\alpha_j$ to $i+1$ and $2n$ back to $1$.

For a sequence $\alpha=(\alpha_1,\ldots,\alpha_p)$ of index sets,
define $\sigma(\alpha)=(\sigma(\alpha_1),\ldots,\sigma(\alpha_p)).$
\end{defn}

\begin{lem}\label{lem:cyclicshift}
Let $A$ be a $n \times n$ totally positive matrix.  Then there exists
 a totally positive $n \times n$ matrix $B$ and a positive constant $c_A$ such that
\[ \p{\sigma(\alpha_j)}(B) = c_A\, \p{\alpha_j}(A) \] for all index sets $\alpha_j$, where $\sigma$ is the cyclic shift operator as defined
in Definition~\ref{defn:cyclicshift}.

In particular, if $\alpha=(\alpha_1,\ldots,\alpha_p)$ and
$\beta=(\beta_1,\ldots,\beta_p)$ then $\alpha/\beta$ is bounded if
and only if $\sigma(\alpha)/\sigma(\beta)$ is bounded.
\end{lem}

\begin{pf}
Let $\bar A$ be the standard matrix representation of the embedding of
$A$ into $TPGr(n,2n)$ as discussed in
Section~\ref{sec:grassmannian}.  Enumerate the rows of $\bar A$ as $\bar A_1, \bar A_2, \ldots, \bar A_{2n}$.

Form the element $\Lambda \in TPGr(n,2n)$ which is represented by the matrix $C$ having rows
\[
	C = [(-1)^{n-1} \bar A_{2n};\, \bar A_1;\, \bar A_2;\, \bar A_3;\, \ldots;\, \bar A_{2n-1}].
\]
Let $\bar B$ be the standard matrix representation of $\Lambda$, i.e. $\bar B = C X$ for some $X\in GL(n,\R)$ where $\det X > 0$.

Then for any index set $\alpha_j$, we have
\begin{eqnarray*}
	[\alpha_j](A) &=& \det \bar A(\alpha_j|(1,2,\ldots,n)) \\
	&=& \det C(\sigma(\alpha_j)|(1,2,\ldots,n)) \\
	&=& \det \bar B(\sigma(\alpha_j)|(1,2,\ldots,n)) \cdot \det X^{-1} \\
	&=& [\sigma(\alpha_j)](B) \cdot \det X^{-1}
\end{eqnarray*}
where $B$ is the totally positive matrix corresponding to $\Lambda$.
\end{pf}

Analogous to the cyclic shift operator is the reversal operator which is described by Fallat et al.~but which behaves differently in this situation of non-principal minors (see \cite[\textsection 3]{Fallat:2003az}).

\begin{defn}[Reversal]
\label{defn:reversal} For an index set $\alpha_j$, define the
\emph{reversal} of the elements of $\alpha_j$ as
\[ \rho(\alpha_j) = \{ (2n+1) -i \, | \, i \in \alpha_j \}. \]

For a sequence $\alpha=(\alpha_1,\ldots,\alpha_p)$ of index sets,
define $\rho(\alpha)=(\rho(\alpha_1),\ldots,\rho(\alpha_p))$
\end{defn}

\begin{lem}
Let $A$ be a $n \times n$ totally positive matrix.  Then there exists
 a totally positive $n \times n$ matrix $B$ and a positive constant $c_A$ such that
\[ \p{\rho(\alpha_j)}(B) = c_A\, \p{\alpha_j}(A) \] for all index sets $\alpha_j$ where $\rho$ is the reversal operator as defined in Definition~\ref{defn:reversal}.

In particular, if $\alpha=(\alpha_1,\ldots,\alpha_p)$ and
$\beta=(\beta_1,\ldots,\beta_p)$ then $\alpha/\beta$ is bounded if
and only if $\rho(\alpha)/\rho(\beta)$ is bounded.
\label{lem:reversal}
\end{lem}

We leave the details of the proof to the reader.

\section{Necessary Conditions for Bounded Ratios}

For $i\in\{1,2,\ldots,2n\}$, let
$f_\alpha(i)$ be the number of index sets in $\alpha$ that contain
$i$.  We now give a generalization of a simple, necessary, but not
sufficient condition for a ratio to be bounded originally described by Fallat et al.~(see
\cite{Fallat:2003az,Skandera:2004qp}).

\begin{defn}[ST0 Condition]
Let $\alpha$ and $\beta$ be two sequences of index sets. If
$f_\alpha(i)=f_\beta(i)$ for all $i$, we say the ratio
$\alpha/\beta$ satisfies the ST0 (set-theoretic) condition.
\end{defn}

\begin{prop}\label{prop:ST0} Let
$\alpha=(\alpha_1,\ldots,\alpha_p)$ and
$\beta=(\beta_1,\ldots,\beta_p)$ be two sequences of index sets,
with each set containing the same number of elements. If $\alpha / \beta$ is
bounded for all totally positive matrices, then the ratio satisfies
the ST0 condition.
\end{prop}

\begin{pf}
Suppose that $\alpha / \beta$ does not satisfy the ST0 condition.  By Lemma~\ref{lem:cyclicshift}, we may assume without loss of generality that $f_\alpha(1) \neq f_\beta(1)$.

Let $C$ be any totally positive matrix, for example the matrix arising from Figure~\ref{fig:general} when all weights are chosen to be $1$.  Let $\bar C$ be the standard matrix representation of the embedding of $C$ into $TPGr(n,2n)$, and enumerate the rows of $\bar C$ as $\bar C_1, \bar C_2, \ldots, \bar C_{2n}$.

Construct a new element in $TPGr(n,2n)$ which has matrix representative $\bar A$ whose rows are
\[
	\bar A = [t \bar C_1;\, \bar C_2;\, \bar C_3;\, \ldots;\, \bar C_{2n}]
\]
where $t$ is chosen to be a positive indeterminate.  (We can think of $\bar A$ as the embedding of $A = \text{diag}(t,1,1,\ldots,1)\times C$ into $TPGr(n,2n)$).

Let $\alpha_j$ be any index set.  Then either:
\begin{itemize}
	\item $1\in \alpha_j$ and $[\alpha_j](A) = c_i t$ for some positive constant $c_i$; or
	\item $1\not\in \alpha_j$ and $[\alpha_j](A) = c_i$ for some positive constant $c_i$.
\end{itemize}

Thus $\alpha(A)$ is a monomial in $t$ of degree $f_\alpha(1)$ and $\beta(A)$ is a monomial in $t$ of degree $f_\beta(1)$.  Because we have assumed that $f_\alpha(1) \neq f_\beta(1)$, $\alpha(A)/\beta(A)$ must increase without bound as either $t\to 0$ or $t\to\infty$ so the ratio is not bounded. 
\end{pf}

In order to present another necessary condition for a ratio $\alpha / \beta$ to be bounded, we first need to discuss the concept of
majorization. The following two definitions and one proposition can be
found in \cite{Marshall:1979az}.

\begin{defn}[Majorization]
Let $x=(x_1,\ldots,x_n)$ and $y=(y_1,\ldots,y_n)$ be two
non-increasing sequences of non-negative integers.  Then $x$
majorizes $y$ (written $x\succeq y$) if for each $k=1,2,\ldots,n$
\[ \sum_{i=1}^k x_i \geq \sum_{i=1}^k y_i \] with equality if $k=n$.
\end{defn}

\begin{defn}[Conjugate Sequence]
The conjugate sequence to $x=(x_1,\ldots,x_n)$ is given by
$x^*=(x^*_1,\ldots,x^*_{x_1})$ where
\[ x^*_j = |\{i : x_i\geq j \}|. \]
\end{defn}

\begin{prop}
Let $x=(x_1,\ldots,x_n)$ and $y=(y_1,\ldots,y_n)$ be two
non-increasing sequences of non-negative integers.  Then $x\succeq
y$ if and only if $y^* \succeq x^*$.
\end{prop}

The following notion of an interval is relied upon in the work of Fallat, Skandera, and others (see \cite{Fallat:2003az,Skandera:2004qp}).  Note that we define an interval slightly differently as contiguous points on a labeled $2n$-gon rather than contiguous points on a line segment with $2n$ vertices, but Lemma~\ref{lem:cyclicshift} shows us that such a distinction is irrelevant in most cases.

\begin{defn}[Interval]
A subset $L\subseteq \{1,2,\ldots,2n\}$ is called an interval if
either $L$ or $L^c=\{1,2,\ldots,2n\}\setminus L$ has the form
$\{i,i+1,i+2,\ldots,i+m\}$.
\end{defn}

Unless mentioned otherwise, all intervals $L$ will be assumed to be
subsets of $\{1,\ldots,2n\}$ of the specified form.

The following definition comes directly from the work of Fallat et al.~but is applied to the case of non-principal minors (see \cite[\textsection 2]{Fallat:2003az}).

\begin{defn}[Condition (M)]\label{def:m}
Let $\alpha=(\alpha_1,\alpha_2,\ldots,\alpha_p)$ and
$\beta=(\beta_1,\beta_2,\ldots,\beta_p)$ be two sequences of index
sets. For any subset $L$ of $\{1,\ldots,2n\}$, define $m(\alpha,L)$
to be the non-increasing rearrangement of the sequence
$(|\alpha_1\cap L|,\ldots, |\alpha_p\cap L|)$.  We say that a ratio
$\alpha / \beta$ satisfies condition (M) if
\[ m(\alpha,L) \succeq m(\beta,L) \]
for every interval $L$.
\end{defn}

\begin{rem}
	Condition (M) implies the ST0 condition by choosing the interval $L=\{j\}$ for $j=1,\dotsc,2n$.
\end{rem}

Before we show that condition (M) is necessary for a ratio $\alpha/\beta$ to be bounded, we give some lemmas which will aid in the proof.

\begin{lem}\label{lem:L<=>L^c}
Let $\alpha=(\alpha_1,\alpha_2,\ldots,\alpha_p)$ and
$\beta=(\beta_1,\beta_2,\ldots,\beta_p)$ be two sequences of index
sets and $L$ be any interval.  If $m(\alpha,L) \succeq m(\beta,L)$
then $m(\alpha,L^c) \succeq m(\beta,L^c)$.
\label{mcomplement}\end{lem}

\begin{pf}
Denote the components of $m(\alpha,L)$, etc. by
\begin{eqnarray*}
m(\alpha,L)&=&(m_1,\ldots,m_p), \\
m(\beta,L)&=&(m'_1,\ldots,m'_p), \\
m(\alpha,L^c)&=&(n_p,\ldots,n_1),\\
m(\beta,L^c)&=&(n'_p,\ldots,n'_1). \\
\end{eqnarray*}

Since $|\alpha_j\cap L|+|\alpha_j\cap L^c|=n$, we have
$m_i+n_i=m'_i+n'_i=n$ for all $i$.  Let $M = \sum_{i=1}^p m_i = \sum_{i=1}^p m'_i$.  Then for any index $k\leq p$ we
have
\begin{eqnarray*}
n_p+n_{p-1}+\cdots+n_{p-(k-1)} &=& (n-m_p)+\cdots+(n-m_{p-(k-1)}) \\
&=& k n-M + m_1+m_2+\cdots+m_{p-k} \\
&\geq& kn - M + m'_1+m'_2+\cdots+m'_{p-k} \\
&\geq& n'_p+n'_{p-1}+\cdots+n'_{p-(k-1)}
\end{eqnarray*}
and hence $m(\alpha,L^c) \succeq m(\beta,L^c).$
\end{pf}

Thus a ratio $\alpha / \beta$ satisfies condition $(M)$ if and only
if $m(\alpha,L)\succeq m(\beta,L)$ for all intervals $L$ with
$|L|\leq n$.

The following theorem is a direct analog of a theorem of Fallat et al.~and is proved in a similar fashion (see \cite[Theorem 2.4]{Fallat:2003az}).

\begin{thm}\label{thm:condm}
Let $\alpha=(\alpha_1,\alpha_2,\ldots,\alpha_p)$ and
$\beta=(\beta_1,\beta_2,\ldots,\beta_p)$ be two sequences of index
sets. If the ratio $\alpha / \beta$ is bounded for all totally
positive matrices, then it satisfies condition (M).
\end{thm}

\begin{pf}

Let $\alpha = (\alpha_1,\alpha_2,\ldots,\alpha_p)$ and $\beta =
(\beta_1,\beta_2,\ldots,\beta_p)$ be two sequences of index sets
such that $\alpha / \beta$ is bounded. 
By Lemmas~\ref{lem:cyclicshift} and \ref{lem:L<=>L^c}, it is sufficient to show that $m(\alpha,L) \succeq m(\beta,L)$ for all intervals $L=\{1,\ldots,s\}$ with $s\leq n$.

Fix $s \leq n$ and let $L=\{1,\ldots,s\}$.  We then construct totally
positive matrices $A_1,A_2,\ldots,A_s$ as follows:

Consider the
planar network shown in Figure~\ref{majorization
  diagram}.  Define the matrix $A_i$ to be the matrix corresponding to
this planar network (see Section~\ref{subsec:planar_networks}) with weights $a_1=a_2=\cdots=a_i=t$ where $t$ is a positive indeterminate and remaining weights $a_{i+1}=\cdots=a_s=1$.

\begin{figure}[h!]
\centering
\includegraphics[width=5.4in]{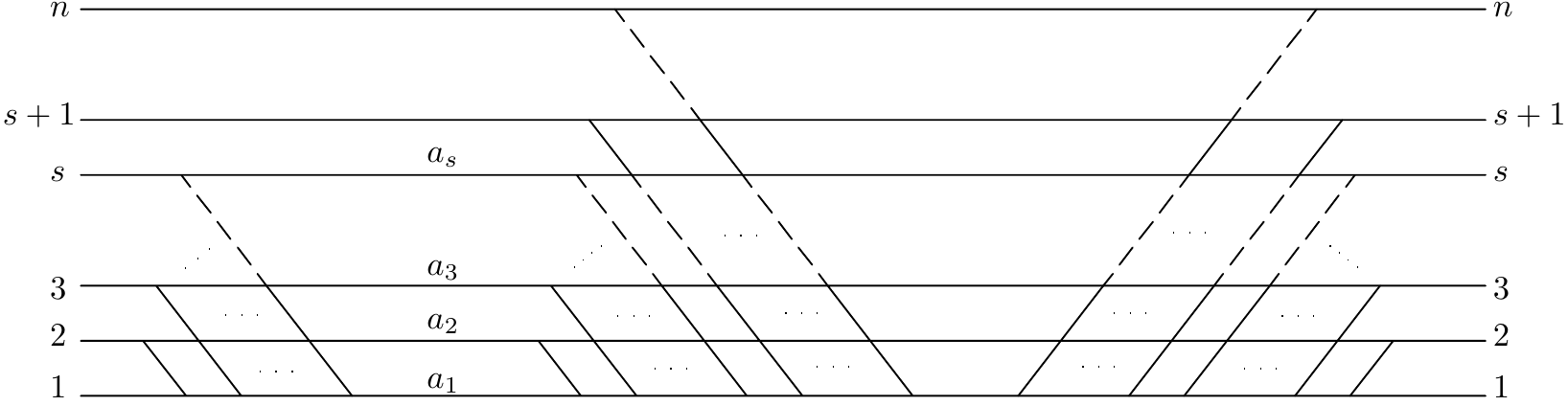}
\caption{Diagram for Matrices $A_i$}\label{majorization diagram}
\end{figure}

If $\alpha_i$ is any index set, then $[\alpha_i](A_k)$ is a polynomial in $t$ and has a well-defined degree.
In fact, recalling that $[\alpha_i](A_k)$ is the sum of weights of path families with sources at $\alpha_i\cap\{1,\ldots,n\}$ and 
sinks at $\{n+1,\ldots,2n\}\setminus \alpha_i$ (see Section~\ref{subsec:planar_networks}), we have $\deg[\alpha_i](A_k) = \min(k, |\alpha_i\cap L |)$
and hence
\[
	\deg \alpha(A_k) = \sum_{i=1}^p \deg [\alpha_i](A_k) = \sum_{i=1}^p \min(k,|\alpha_i\cap L |).
\]

For $\alpha / \beta$ to be bounded as $t\to\infty$, we must have that $\deg
\alpha(A_k)\leq \deg \beta(A_k)$, i.e.
\[\sum_{i=1}^p \min(k, |\alpha_i\cap L)|)
\leq
\sum_{i=1}^p \min(k,|\beta_i\cap L|)\]
for each $1\leq k \leq s$.

Note that if $m^*(\alpha,L) = (m_1^*(\alpha,L), m_2^*(\alpha,L), \ldots)$ is the conjugate sequence to $m(\alpha,L)$ we recognize the left side of the inequality as
\[
	\sum_{i=1}^p \min(k, |\alpha_i\cap L)|) = \sum_{j=1}^k \left|\left\{ i : |\alpha_i \cap L| \geq j \right\}\right| = \sum_{j=1}^k m_j^*(\alpha,L).
\]
and similarly for the summation with $\beta$.

Since this inequality holds for all $1\leq k\leq s$ with equality when $k=s$ by
the ST0 condition (see Proposition~\ref{prop:ST0}), we have that $m^*(\alpha,L) \preceq m^*(\beta,L)$
and hence $m(\alpha,L) \succeq m(\beta,L)$ as desired.
\end{pf}

\section{Basic and Elementary Bounded Ratios}

In this section we define two special classes of ratios of the form \[
\plucker{\alpha_1}{\alpha_2}{\beta_1}{\beta_2}
\] where the $\alpha_1$, $\alpha_2$, $\beta_1$, and $\beta_2$ are index sets.  In particular, ratios belonging to either of these classes will be bounded by 1.

Let $\alpha_i$ be an index set.  If $\alpha_i=\gamma_1\cup\cdots
\cup \gamma_n$ then the notation $[\gamma_1,\ldots,\gamma_n]$ should
be interpreted as $[\alpha_i]$.  Furthermore, if $\gamma_i$ consists of only
a single element $j$, then we may simply write $j$ instead of
$\gamma_i$.

\begin{defn}
An elementary ratio is a ratio of the form
\[\plucker{i,j',\Delta}{i',j,\Delta}{i,j,\Delta}{i',j',\Delta}
\] satisfying 
\begin{enumerate}
	\item $|\Delta|=n-2$;
	\item $i<i'<j<j'$ when considering each element as the mod $2n$ representative in $\{i,i+1,\ldots,i+2n-1\}$; and
	\item $i$, $i'$, $j$, and $j'$ are not elements of $\Delta$.
\end{enumerate}
\end{defn}

\begin{prop}
	A ratio $R$ of the form 
	\[
		\plucker{i,j',\Delta}{i',j,\Delta}{i,j,\Delta}{i',j',\Delta}
	\]
	with $|\Delta|=n-2$, $\Delta\cap \{i,i',j,j'\}=\emptyset$ and $i$, $i'$, $j$, $j'$ pairwise distinct is elementary if and only if it satisfies condition (M).
	\label{prop:elementary<=>conditionM}
\end{prop}

\begin{pf}
	Suppose the ratio $R$ is elementary, and let $L$ be any interval.  
	Set $\alpha_1=\{i,j'\}\cup \Delta$, $\alpha_2 = \{i',j\}\cup \Delta$, $\beta_1 = \{i,j\}\cup \Delta$, and $\beta_2 = \{i',j'\}\cup \Delta$, so that $R=[\alpha_1][\alpha_2]/[\beta_1][\beta_2]$.
	
	Because $R$ satisfies the ST0 condition by construction, it suffices to verify that 
	\[
		\max\left( |\alpha_1 \cap L|, |\alpha_2 \cap L| \right) \geq \max\left( \beta_1 \cap L|,|\beta_2 \cap L|\right)
	\]
	or equivalently
	\[
		\max\left( |\{i,j'\} \cap L|, |\{i',j\} \cap L| \right) \geq \max\left( \{i,j\} \cap L|,| \{i',j'\} \cap L|\right)
	\]
	noting that $|\Delta \cap L|$ appears in every term and thus may be omitted.
	We verify this last inequality by considering the possible values for the right hand side.
	
	If the right hand side is $0$, the inequality is trivially satisfied.
	
	If the right hand side is $1$, the interval $L$ contains at least one of $i$, $i'$, $j$, or $j'$ and thus the left hand side is at least $1$.
	
	If the right hand side is $2$, the interval $L$ contains either $i$ and $j$, or $i'$ and $j'$.  Assume for the moment that $L$ contains both $i$ and $j$.  Then because $R$ is an elementary ratio it must be that the interval $L$ also contains either $i'$ or $j'$ (or possibly both), and hence the left hand side is $2$.  Similar reasoning holds if $L$ had instead contained $i'$ and $j'$.

	Conversely, suppose that the ratio $R$ satisfies condition (M).  Consider the two intervals $L = \{i, i+1, \dotsc, j\}$ and $L' = \{j, j+1, \dotsc, i\}$, working with the elements modulo $2n$ as required.  Because $R$ satisfies condition (M), it must be that 
\[
	\max\left( |\{i,j'\} \cap L|, |\{i',j\} \cap L| \right) \geq \max\left( \{i,j\} \cap L|,| \{i',j'\} \cap L|\right) = 2
\]
	and hence either $i'$ or $j'$ lies in $L$.  Additionally, upon consideration of the complementary interval $L'$, we see that 
\[
	\max\left( |\{i,j'\} \cap L|, |\{i',j\} \cap L| \right) \geq \max\left( \{i,j\} \cap L|,| \{i',j'\} \cap L|\right) = 2
\]
	and hence either $i'$ or $j'$ lies in $L'$.  Thus working modulo $2n$ and considering representatives in $\{i, i+1, \dotsc , i+2n-1\}$, we have either $i < i' < j < j'$ or $i < j' < j < i'$.  In the first case the ratio is elementary.  In the latter case, a simple renaming $i \leftrightarrow j'$ and $i' \leftrightarrow j$ preserves the ratio and makes it elementary.
\end{pf}


\begin{rem}
All elementary ratios are necessarily bounded by 1. Indeed, the
short Pl\"ucker relation
\begin{equation}\p{i,i',\Delta}\p{j,j',\Delta}+\p{i,j',\Delta}\p{i',j,\Delta}=\p{i,j,\Delta}\p{i',j',\Delta}\label{shortplucker}\end{equation}
together with the positivity of all Pl\"ucker coordinates over $TPGr(n,2n)$ imply
\[\p{i,j',\Delta}\p{i',j,\Delta}<\p{i,j,\Delta}\p{i',j',\Delta}\] as desired.  
\end{rem}

Computationally the elementary ratios are inefficient due to the large number of them. The solution to this problem is to consider instead a small subset of the elementary ratios, which we will call the {\em basic} ratios.  We will show that every
elementary ratio can be written as a product of positive powers of
basic ratios.  We will use this fact in the next section.

\begin{defn}A basic ratio is one of the form
\[ \plucker{i,j+1,\Delta}{i+1,j,\Delta}{i,j,\Delta}{i+1,j+1,\Delta} \]
where $i,j\in \{1,\ldots,2n\}$ and $\Delta\subset \{1,\ldots,2n\}$
such that $|\Delta|=n-2$ and $i$, $i+1$, $j$, $j+1$ and $\Delta$ are
all distinct.  Here indices $i+1$ and $j+1$ are understood mod $2n$.
\end{defn}

Clearly, a basic ratio is an elementary ratio with $i'=i+1$, and $j'=j+1$.  

We define the complexity of a particular elementary ratio 
\[ R = \plucker{i,j',\Delta}{i',j,\Delta}{i,j,\Delta}{i',j',\Delta} \]
 as
\[ 
	\mu (R) = \left| \Delta \cap
	\left(\{i,i+1,\ldots,i'\} \cup \{j,j+1,\ldots,j'\}\right) \right|.
\]
To prove that every
elementary ratio can be written as a product of basic ratios, we first
consider the following special case.
\begin{lem}\label{lem:mu=0}
	An elementary ratio $R$ with complexity $\mu(R)=0$ can be written as a product of basic ratios.
\end{lem}
\begin{pf}
	We define 
	\[
		\delta(R) = \left| \{i,i+1,\ldots,i'\} \cup \{j,j+1,\ldots,j'\} \right|.
	\]
	Recall that $R$ is a basic ratio if $i' = i+1$ and $j' = j+1$, or in other words $\delta(R) = 4$.  We proceed by induction.

	Assume that when $\mu(R) = 0$ and $\delta(R) = 4,5,\ldots,k-1$ we have a factorization of ratio $R$ into a product of basic ratios.
	
	Now consider a given elementary ratio $R$ with $\mu(R) = 0$ and $\delta(R) = k > 4$.
	It cannot be the case that both $i' = i+1$ and $j' = j+1$ as $\delta(R)>4$.  Without loss of generality, assume that $i+1 \neq i'$ (otherwise exchange the labels of $i$ and $i'$ with $j$ and $j'$ respectively).
	
	Now the elementary ratio $R$ factors as
	\[ \plucker{i,j',\Delta}{i',j,\Delta}{i,j,\Delta}{i',j',\Delta} = \left(\plucker{i,j',\Delta}{i+1,j,\Delta}{i,j,\Delta}{i+1,j',\Delta}\right)\left(\plucker{i+1,j',\Delta}{i',j,\Delta}{i+1,j,\Delta}{i',j',\Delta}\right) \]
	where each factor $F$ on the right hand size has $\mu(F)=0$ and $\delta(F) < k$.  By induction, each factor on the right hand side can be expressed as a product of basic ratios.  Hence $R$ can we written as a product of basic ratios.
\end{pf}

\begin{thm}
Every elementary ratio can be written as a product of basic
ratios.\label{jawtheorem}
\end{thm}

\begin{pf}
We shall proceed by induction on $\mu(R)$.  	
By the previous lemma, when $\mu(R) = 0$ we have that $R$ can be expressed as a product of basic ratios.
Assume that for any ratio $R$ with $\mu(R) = 0,1,\ldots,k-1$ we can express $R$ as a product of basic ratios.

Now consider an elementary ratio $R$ with $\mu(R) = k > 0$.  It cannot be the case that both $\Delta \cap \{i,i+1,\ldots,i'\} = \emptyset$ and $\Delta \cap \{j,j+1,\ldots,j'\} = \emptyset$, so assume without loss of generality that $\Delta \cap \{i,i+1,\ldots,i'\} \neq \emptyset$ (if not, exchange the labels of $i$ and $i'$ with $j$ and $j'$ respectively).

Let $p \in \Delta \cap \{i,i+1,\ldots,i'\}$ be the element
nearest to $i$ and let $\Delta'=\Delta\setminus\{p\}$.  

The ratio $R$ then factors as
\[ \plucker{i,j',\Delta}{i',j,\Delta}{i,j,\Delta}{i',j',\Delta} = \left(\plucker{i,j',p,\Delta'}{i,i',j,\Delta'}{i,j,p,\Delta'}{i,i',j',\Delta'}\right)\left(\plucker{i,i',j',\Delta'}{i',j,p,\Delta'}{i,i',j,\Delta'}{i',j',p,\Delta'}\right)
\]
where each factor $F$ on the right hand side has $\mu(F) < k$ and hence may be written as a product of basic ratios.
\end{pf}

\section{A Factorization of (Some) Bounded Ratios}

In this section we give an alternative proof of a necessary and
sufficient condition for a ratio of the form \begin{equation}
\plucker{\alpha_1}{\alpha_2}{\beta_1}{\beta_2}\label{eq:2over2}
\end{equation} to be bounded in terms of the four index sets,
$\alpha_1$, $\alpha_2$, $\beta_1$, and $\beta_2$.  In addition, we
will show that this ratio can be written as a product of elementary ratios.

For the remainder of the section we will assume that $R$ is a ratio of the form $\iplucker{\alpha_1}{\alpha_2}{\beta_1}{\beta_2}$ which satisfies the ST0 condition and condition (M).  Denote the set of all such ratios by $\mathcal B$.
We define
\begin{align*}
\Delta &= \alpha_1\cap\alpha_2=\beta_1\cap\beta_2;   \\
\gamma_1 &= (\alpha_1 \cap \beta_1)\setminus\Delta;  \\
\gamma_2 &= (\alpha_1 \cap \beta_2)\setminus\Delta;  \\
\delta_1 &= (\alpha_2 \cap \beta_2)\setminus\Delta;  \\
\delta_2 &= (\alpha_2 \cap \beta_1)\setminus\Delta;  \text{ and}  \\
\Omega &= \gamma_1 \cup \gamma_2 \cup \delta_1 \cup \delta_2, 
\end{align*}
so that
\[ R \; = \; \plucker{\alpha_1}{\alpha_2}{\beta_1}{\beta_2} \; = \; \plucker{\gamma_1,\gamma_2,\Delta}{\delta_1,\delta_2,\Delta}{\gamma_1,\delta_2,\Delta}{\delta_1,\gamma_2,\Delta}. \]
(Recall that notationally $[\gamma_1,\gamma_2,\Delta]$ means $[\gamma_1\cup\gamma_2\cup\Delta]$, and that we necessarily have: $|\gamma_1| = |\delta_1|$;
$|\gamma_2| = |\delta_2|$; and  $|\gamma_1|+|\gamma_2|+|\Delta| =
n$.)

An important property of the ratio $R$ is the number of indices which are \emph{not} shared by all index sets comprising the ratio.  We shall denote this quantity as 
\begin{align*}
	\nu(R) &= n - |\Delta| = |\Omega|/2.
\end{align*}

Before proceeding, we investigate what information $\nu(R)$ holds.

\begin{defn}[Trivial Ratio]
	We say a ratio $\iplucker{\alpha_1}{\alpha_2}{\beta_1}{\beta_2}$ is trivial if either
	\begin{itemize}
		\item $\alpha_1 = \beta_1$ and $\alpha_2 = \beta_2$; or
		\item $\alpha_1 = \beta_2$ and $\alpha_2 = \beta_1$.
	\end{itemize}
\end{defn}

Note that a ratio $R$ satisfying the ST0 condition with $\nu(R) = 0$ or $1$ is trivial.

\begin{lem}\label{lem:nu(R)=2=>trivial-or-elementary}
	Suppose $ R\in \mathcal B$ and $\nu(R) = 2$.  Then either 
	\begin{itemize}
		\item $R$ is trivial; or 
		\item $R$ is an elementary ratio and can be written as a product of basic ratios.
	\end{itemize}\label{nu=2 easy}
\end{lem}

\begin{pf}
	If $R$ is not trivial, $R$ must be of the form 
	\[ \plucker{i,j',\Delta}{i',j,\Delta}{i,j,\Delta}{i',j',\Delta} \]
	with $|\Delta| = n-2$ and $i$, $i'$, $j$, $j'$, and $\Delta$ pairwise distinct.  Proposition~\ref{prop:elementary<=>conditionM} establishes that $R$ is an elementary ratio, and hence $R$ may be written as a product of basic ratios by Theorem~\ref{jawtheorem}.
\end{pf}

We will eventually show that any ratio $R\in \mathcal B$ with $\nu(R) \ge 3$
can be factored
as $R=R_1 R_2$ with $R_i \in \mathcal B$ and $\nu(R_i) < \nu (R)$ for $i=1,2$.  In order to do this, we will rely heavily upon the following definition, simple remark, and technical lemma.

\begin{defn}[Interlacing]Suppose $(i_1,i_2,\ldots,i_k)$ and $(j_1,j_2,\ldots,j_k)$ are two
subsequences of $(1,2,\ldots,2n)$.  Then we say the sequence
$(i_s)$ interlaces the sequence $(j_t)$ if either:
\begin{enumerate}
\item $j_1\leq i_1\leq j_2\leq i_2\leq\cdots\leq j_k\leq i_k$; or
\item $i_1\leq j_1\leq i_2\leq j_2\leq\cdots\leq i_k\leq j_k$.
\end{enumerate}
\end{defn}

\begin{rem}\label{rem:interlacing_denominator_implies_M}
	Suppose the ratio $R = \iplucker{\alpha_1}{\alpha_2}{\beta_1}{\beta_2}$ satisfies the ST0 condition, and suppose that $\beta_1 \setminus \Delta$ interlaces with $\beta_2 \setminus \Delta$.  Then $R$ automatically satisfies condition (M).
\end{rem}

For notational convenience, let $g(\alpha_1,\alpha_2,L) =
\max(|\alpha_1 \cap L|, |\alpha_2 \cap L|)$.  

\begin{lem}[Technical Lemma] \label{lem:technical_lemma}
Suppose that
$ R\in \mathcal B$ and we have a factorization of $R$ as
\begin{equation*}
\plucker{\gamma_1,\gamma_2,\Delta}{\delta_1,\delta_2,\Delta}{\gamma_1,\delta_2,\Delta}{\delta_1,\gamma_2,\Delta}
=
\plucker{\gamma_1,\gamma_2,\Delta}{\gamma_{11},\delta_{12},\delta_2,\Delta}{\gamma_1,\delta_2,\Delta}{\gamma_{11},\delta_{12},\gamma_2,\Delta}
\cdot
\plucker{\gamma_{11},\delta_{12},\gamma_2,\Delta}{\delta_1,\delta_2,\Delta}{\gamma_{11},\delta_{12},\delta_2,\Delta}{\delta_1,\gamma_2,\Delta}
\end{equation*}
for some non-empty sets $\gamma_{11}$, $\gamma_{12}$, $\delta_{11}$, and $\delta_{12}$ where $\gamma_1=\gamma_{11}\cup\gamma_{12}$ and
$\delta_1=\delta_{11}\cup\delta_{12}$ such that $|\gamma_{11}|=|\delta_{11}|$ and $|\gamma_{12}|=|\delta_{12}|$.

Suppose as well that $\gamma_{11} \cup \delta_{12}$ interlaces with $\gamma_{12} \cup \delta_{11}$.

 Then each of the factors
\[R_1 = \plucker{\gamma_1,\gamma_2,\Delta}{\gamma_{11},\delta_{12},\delta_2,\Delta}{\gamma_1,\delta_2,\Delta}{\gamma_{11},\delta_{12},\gamma_2,\Delta}
~~\text{and}~~
R_2 = \plucker{\gamma_{11},\delta_{12},\gamma_2,\Delta}{\delta_1,\delta_2,\Delta}{\gamma_{11},\delta_{12},\delta_2,\Delta}{\delta_1,\gamma_2,\Delta},\]
are elements of $\mathcal B$, and $\nu(R_i) < \nu(R)$ for $i=1,2$.

\label{lem22}
\end{lem}

\begin{pf}  
Observe that
\[
	|\gamma_1 \cap L| + |\delta_1 \cap L| = |(\gamma_{11} \cup \delta_{12}) \cap L| + |(\gamma_{12} \cup \delta_{11}) \cap L|
\] for all intervals $L$.  This, along with the hypothesis that $\gamma_{11} \cup \delta_{12}$ interlaces with $\gamma_{12} \cup \delta_{11}$, immediately gives
\[	g(\gamma_1,\delta_1,L) \geq g(\gamma_{11} \cup \delta_{12},\gamma_{12} \cup \delta_{11},L)
\]	for all intervals $L$.

Fix an interval $L$. Then there are three possible cases:
\begin{enumerate}
\item $|\alpha_1\cap L|>|\alpha_2\cap L|$;
\item $|\alpha_1\cap L|<| \alpha_2\cap L|$; or
\item $|\alpha_1\cap L|=| \alpha_2\cap L|$.
\end{enumerate}
Suppose case (1) holds.  Then since $g(\alpha_1,\alpha_2, L)\geq g(\beta_1,\beta_2,L)$ it follows that $|\gamma_1\cap L|\geq
|\delta_1\cap L|$ and $|\gamma_2\cap L|\geq |\delta_2\cap L|$.  However, if $g(\gamma_1,\delta_1,L) \geq g(\gamma_{11} \cup \delta_{12},\gamma_{12} \cup \delta_{11},L)$ and
$|\gamma_1\cap L|\geq |\delta_1\cap L|$, then applying similar reasoning reveals that
$|\gamma_{11}\cap L|\geq |\delta_{11}\cap L|$
and
$|\gamma_{12}\cap L|\geq |\delta_{12}\cap L|$.  Thus the following four inequalities hold:
\begin{enumerate}
\item[(i)] $|\gamma_1\cap L|\geq |\delta_1\cap L|$;
\item[(ii)] $|\gamma_2\cap L|\geq |\delta_2\cap L|$;
\item[(iii)] $|\gamma_{11}\cap L|\geq |\delta_{11}\cap L|$; and
\item[(iv)] $|\gamma_{12}\cap L|\geq |\delta_{12}\cap L|$.
\end{enumerate}
But these preceding inequalities (i)-(iv) imply that both of the ratios
\[\plucker{\gamma_1,\gamma_2,\Delta}{\gamma_{11},\delta_{12},\delta_2,\Delta}{\gamma_1,\delta_2,\Delta}{\gamma_{11},\delta_{12},\gamma_2,\Delta}
~~\text{and}~~
\plucker{\gamma_{11},\delta_{12},\gamma_2,\Delta}{\delta_1,\delta_2,\Delta}{\gamma_{11},\delta_{12},\delta_2,\Delta}{\delta_1,\gamma_2,\Delta},\]
satisfy condition (M), for the fixed interval $L$.  Similar analysis holds for cases (2) and (3) when
$|\alpha_1\cap L|\leq |\alpha_2\cap L|$ and is omitted here.

Lastly, note that $\nu(R_1) = \nu(R) - |\gamma_{11}| < \nu(R)$ and $\nu(R_2) = \nu(R)-|\delta_{12}| < \nu(R)$.
\end{pf}


\begin{lem}\label{lem:some_dont_interlace}
	Let $ R\in \mathcal B$, and suppose that either
	\begin{itemize}
		\item $\gamma_1$ and $\delta_1$ do not interlace; or
		\item $\gamma_2$ and $\delta_2$ do not interlace (or both).
	\end{itemize}
	Then we may write $R=R_1R_2$ for some ratios $ R_1,R_2\in \mathcal B$ with $\nu(R_i)<\nu(R)$ for $i=1,2$.
\end{lem}

\begin{pf}
	Without loss of generality, assume that $\gamma_1$ and $\delta_1$ do not interlace.  (If instead $\gamma_2$ and $\delta_2$ do not interlace, interchange the labeling of $\alpha_1$ and $\alpha_2$).
	
Label the elements of $\gamma_1\cup \delta_1$ as $\{i_1,i_2,\ldots,i_{2m}\}$ with $i_1<i_2<
\cdots<i_{2m}$, and define $\gamma_{11}=\gamma_1\cap\{i_1,i_3,\ldots,i_{2m-1}\}$,
$\gamma_{12} = \gamma_1 \cap \{i_2,i_4,\ldots,i_{2m}\}$, $\delta_{12} = \delta_1
\cap \{i_1,i_3,\ldots,i_{2m-1}\}$, and $\delta_{11} = \delta_1 \cap \{i_2,i_4,\ldots,i_{2m}\}$.

Because $\gamma_1$ does not interlace with $\delta_1$, we have constructed $\gamma_{11}$, $\gamma_{12}$, $\delta_{11}$, and $\delta_{12}$ to all be non-empty.  In addition, $\gamma_{11}\cup\delta_{12}$ interlaces with $\gamma_{12}\cup\delta_{11}$, thus satisfying the requirements of Lemma~\ref{lem:technical_lemma}.
\end{pf}

We now examine the situation when both $\gamma_1$ interlaces with $\delta_1$ and 
$\gamma_2$ interlaces with $\delta_2$.  

\begin{claim}[Agreeable Labeling]\label{claim:wlog_labeling}
	Let $R\in \mathcal B$ be a ratio satisfying condition (M), and suppose that $\gamma_1$ interlaces with $\delta_1$ and $\gamma_2$ interlaces with $\delta_2$.  Set $\Omega=\gamma_1\cup\gamma_2\cup\delta_1\cup\delta_2=\{i_1,i_2,\ldots,i_{2m}\}$ with $i_1<i_2<\cdots <i_{2m}$.
	Then, up to a possible relabeling of $\alpha_1$ and $\alpha_2$ or $\beta_1$ and $\beta_2$, we may assume that 
	\begin{enumerate}
		\item $\beta_1 \setminus \Delta = \gamma_1\cup\delta_2=\{i_1,i_3,\ldots,i_{2m-1}\}$;
		\item $\beta_2 \setminus \Delta = \delta_1\cup\gamma_2=\{i_2,i_4,\ldots,i_{2m}\}$; and
		\item $i_1 \in \gamma_1$.
	\end{enumerate}
	
	We will say that ratio with a labeling satisfying conditions (1)-(3) is agreeably labeled.
\end{claim}

\begin{pf}
	Suppose that either $\gamma_1 \cup \delta_2$ or $\delta_1 \cup \gamma_2$ contained a consecutive pair of elements $i_l,i_{l+1} \in \Omega$.  It cannot be that this pair lies entirely in one of $\gamma_1$, $\gamma_2$, $\delta_1$, or $\delta_2$, as this would violate the interlacing hypotheses.  However if $i_l$ and $i_{l+1}$ lie in different sets, for example $\gamma_1$ and $\delta_2$, consideration of condition (M) with the interval $L=\{i_{l}, \dotsc, i_{l+1}\}$ again leads to a contradiction.  
	
	Thus $\gamma_1 \cup \delta_2$ and $\delta_1 \cup \gamma_2$ contain no consecutive pairs of elements of $\Omega$, and hence we may assume (up to relabeling of $\beta_1$ and $\beta_2$) that both (1) and (2) hold.  Lastly, we may swap the labeling of $\alpha_1$ and $\alpha_2$ if necessary to ensure that $i_1 \in \alpha_1$ and hence (3) holds.
\end{pf}

Under this labeling, the element $i_2$ may be in either $\delta_1$ or $\gamma_2$.  We investigate each case separately.

\begin{lem}\label{lem:i_2-in-delta_1}
	Let $R\in \mathcal B$ with $\nu(R) \geq 3$, and suppose that both $\gamma_1$ interlaces with $\delta_1$ and $\gamma_2$ interlaces with $\delta_2$.  Assume that $R$ is agreeably labeled (see Claim~\ref{claim:wlog_labeling}), and set $\Omega=\gamma_1\cup\gamma_2\cup\delta_1\cup\delta_2=\{i_1,i_2,\ldots,i_{2m}\}$ with $i_1<i_2<\cdots <i_{2m}$.  Furthermore, assume that $i_2 \in \delta_1$.
	
	Then we may write $R=R_1R_2$ for some ratios $R_1,R_2\in \mathcal B$ with $\nu(R_i)<\nu(R)$ for $i=1,2$.	
\end{lem}


\begin{pf}
	Because $R$ is agreeably labeled, we know that $\gamma_1\cup\delta_2=\{i_1,i_3,\ldots\}$ with $i_1 \in \gamma_1$.  Define $k \geq 1$ to be the value so that $\{i_1, i_3, i_5, \dotsc, i_{2k-1} \} \subseteq \gamma_1$ and $i_{2k+1}\in \delta_2$.
Similarly define $l \geq 1$ to be the value so that $\{i_2, i_4, \dotsc, i_{2l}\} \subseteq \delta_1$ and $i_{2l+2} \in \gamma_2$.

To summarize, we have set
\begin{align*}
	\gamma_1 &= \{i_1,i_3,i_5,\ldots,i_{2k-1},\ast\}; \\
	\gamma_2 &= \{i_{2l+2},\ast\}; \\
	\delta_1 &= \{i_2,i_4,\ldots,i_{2l},\ast\}; \text{ and} \\
	\delta_2 &= \{i_{2k+1},\ast\},
\end{align*}
where the use of $\ast$ is understood to represent the remaining elements and is not the same in each instance.

Because of the interlacing hypothesis, we must have either
\begin{enumerate}
\item[(a)] $l=k$; or
\item[(b)] $l=k-1$.
\end{enumerate}

For case (a), let $\gamma_{11}=\{i_1\}$, $\gamma_{12} = \gamma_1 \setminus \{i_1\} = \{i_3,i_5,\ldots,i_{2k-1},\ast\}$, $\delta_{11} = \{i_2\}$ and $\delta_{12} = \delta_1 \setminus \{i_2\} = \{i_4,\ldots,i_{2l},\ast\}$.

We claim $\gamma_{12}\cup\delta_2$ and $\delta_{12}\cup\gamma_2$ interlace, as 
\begin{align*}
	\gamma_{12}\cup\delta_2 &= (\gamma_1 \cup \delta_2) \setminus \{i_1\} = \{i_3,i_5,\dotsc,i_{2m-1} \} \text{ and} \\
	\delta_{12}\cup\gamma_2 &= (\delta_1 \cup \gamma_2) \setminus \{i_2\} = \{i_4,i_6,\dotsc,i_{2m}\}.
\end{align*}

Similarly $\gamma_{11}\cup\delta_{2}$ and $\delta_{11}\cup\gamma_2$ interlace, as 
\begin{align*}
	\gamma_{11}\cup\delta_2 &= \{i_1\} \cup \delta_2 = \{i_1, i_{2k+1}, \ast \} \text{ and} \\
	\delta_{11}\cup\gamma_2 &= \{i_2\} \cup \gamma_2 = \{i_2, i_{2l+2}, \ast \},
\end{align*}
noting $l=k$ and $\gamma_2$ interlaces with $\delta_2$ by hypothesis.

Therefore, we may write
\[
\plucker{\gamma_1,\gamma_2,\Delta}{\delta_1,\delta_2,\Delta}{\gamma_1,\delta_2,\Delta}{\delta_1,\gamma_2,\Delta}
=
\plucker{\gamma_1,\gamma_2,\Delta}{\gamma_{11},\delta_{12},\delta_2,\Delta}{\gamma_1,\delta_2,\Delta}{\gamma_{11},\delta_{12},\gamma_2,\Delta}
\cdot
\plucker{\gamma_{11},\delta_{12},\gamma_2,\Delta}{\delta_1,\delta_2,\Delta}{\gamma_{11},\delta_{12},\delta_2,\Delta}{\delta_1,\gamma_2,\Delta}
\]
where both ratios on the right hand side satisfy condition (M) by Remark~\ref{rem:interlacing_denominator_implies_M}.  Labeling the ratios on the right hand side as $R_1$ and $R_2$ respectively, we see that
$\nu(R_1) = \nu(R) - |\gamma_{11}|$ and $\nu(R_2) = \nu(R) - |\gamma_{12}|$.  Now $|\gamma_{11}| = 1$, so $\nu(R_1) < \nu(R)$.  

If $|\gamma_{12}| \geq 1$ we are finished with this case.  If instead $|\gamma_{12}| = 0$ we deduce that $\gamma_1 = \{i_1\}$, $\gamma_2 = \{i_3, i_5, \dotsc, i_{2m-1}\}$, $\delta_1 = \{i_2\}$, and $\delta_2 = \{i_4, i_6, \dotsc, i_{2m} \}$.  We can then factor $R$ as 
\begin{align*}
\plucker{i_1,i_4,i_6,\ldots,\Delta}{i_2,i_3,i_5,\ldots,\Delta}{i_1,i_3,i_5,\ldots,\Delta}{i_2,i_4,i_6,\ldots,\Delta} = &
\plucker{i_1,i_4,i_6,\ldots,\Delta}{i_2,i_4,i_5,i_7,\ldots,\Delta}{i_2,i_4,i_6,\ldots,\Delta}{i_1,i_4,i_5,i_7,\ldots,\Delta} \\
& \cdot \plucker{i_1,i_4,i_5,i_7,\ldots,\Delta}{i_2,i_3,i_5,\ldots,\Delta}{i_2,i_4,i_5,i_7,\ldots,\Delta}{i_1,i_3,i_5,\ldots,\Delta},
\end{align*}
where the ellipses indicate the the sequence continues with the same parity subscripts.  Note that $\{i_1, i_5, i_7, \dotsc\}$ interlaces with $\{i_2, i_6, i_8, \dotsc\}$ and that $\{i_1, i_3\}$ interlaces with $\{i_2, i_4\}$, hence both ratios on the right hand side satisfy condition (M) by Remark~\ref{rem:interlacing_denominator_implies_M}.  Labeling the ratios on the right hand side as $R_1$ and $R_2$ respectively, we see that $\nu(R_1) = \nu(R) - 1 < \nu(R)$ and $\nu(R_2) = 2 < \nu(R)$.

Now we return to case (b), where $l=k-1$.
Let $\gamma_{11} = \{i_3,i_5,\ldots,i_{2k-1}\}$,
$\gamma_{12} = \gamma_1 \setminus \gamma_{11}$,
$\delta_{11} = \{i_2,i_4,\ldots,i_{2l}\}$, and 
$\delta_{12} = \delta_1 \setminus \delta_{11}$.

We claim $\gamma_{12}\cup\delta_2$ and $\delta_{12}\cup\gamma_2$ interlace, as
\begin{align*}
	\gamma_{12}\cup\delta_2 &= (\gamma_1 \cup \delta_2) \setminus \{i_2, i_3, i_4, \dotsc i_{2l+1} \} \text{ and} \\
	\delta_{12}\cup\gamma_2 &= (\delta_1 \cup \gamma_2) \setminus \{i_2, i_3, i_4, \dotsc i_{2l+1} \},
\end{align*}
noting $\gamma_1 \cup \delta_2$ interlaces with $\delta_1 \cup \gamma_2$ and we have removed a section of consecutive elements of $\Omega$.

Similarly, we claim $\gamma_{11}\cup\delta_2$ and $\delta_{11}\cup\gamma_2$ interlace, as 
\begin{align*}
	\gamma_{11}\cup\delta_2 &= \{i_3,i_5,\ldots,i_{2k-1},i_{2k+1},\ast\} \text{ and} \\
	\delta_{11}\cup\gamma_2 &= \{i_2,i_4,\ldots,i_{2l},i_{2l+2},\ast\},
\end{align*}
noting $l=k-1$ and $\gamma_2$ interlaces with $\delta_2$ by hypothesis.

Therefore, we may write
\[
\plucker{\gamma_1,\gamma_2,\Delta}{\delta_1,\delta_2,\Delta}{\gamma_1,\delta_2,\Delta}{\delta_1,\gamma_2,\Delta}
=
\plucker{\gamma_1,\gamma_2,\Delta}{\gamma_{11},\delta_{12},\delta_2,\Delta}{\gamma_1,\delta_2,\Delta}{\gamma_{11},\delta_{12},\gamma_2,\Delta}
\cdot
\plucker{\gamma_{11},\delta_{12},\gamma_2,\Delta}{\delta_1,\delta_2,\Delta}{\gamma_{11},\delta_{12},\delta_2,\Delta}{\delta_1,\gamma_2,\Delta}
\]
where both ratios on the right hand side satisfy condition (M) by Remark~\ref{rem:interlacing_denominator_implies_M}.  Labeling the ratios on the right hand side as $R_1$ and $R_2$ respectively, we see that $\nu(R_1) = \nu(R) - |\gamma_{11}|$ and $\nu(R_2) = \nu(R) - |\gamma_{12}|$. 

Observe that $i_2 \in \delta_{11}$ so $|\delta_{11}| = |\gamma_{11}| \geq 1$ and hence $\nu(R_1) < \nu(R)$.  Similarly, $i_1 \in \gamma_{12}$ so $|\gamma_{12}| \geq 1$ and hence $\nu(R_2) < \nu(R)$.

This concludes the proof, as we have successfully dealt with both cases (a) and (b).
\end{pf}

\begin{lem}\label{lem:i_2-in-gamma_2}
	Let $ R\in \mathcal B$ with $\nu(R) \geq 3$, and suppose that both $\gamma_1$ interlaces with $\delta_1$ and $\gamma_2$ interlaces with $\delta_2$.  Assume that $R$ is agreeably labeled (see Claim~\ref{claim:wlog_labeling}), and set $\Omega=\gamma_1\cup\gamma_2\cup\delta_1\cup\delta_2=\{i_1,i_2,\ldots,i_{2m}\}$ with $i_1<i_2<\cdots <i_{2m}$.  Furthermore, assume that $i_2 \in \gamma_2$.
	
	Then we may write $R=R_1R_2$ for some ratios $ R_1,R_2 \in \mathcal B$ with $\nu(R_i)<\nu(R)$ for $i=1,2$.
\end{lem}
\begin{pf} First, note that $i_3\in \delta_2$, since otherwise $\{i_1,i_3\}\subset \gamma_1$ and hence
$\gamma_1$ and $\delta_1$ would not interlace.  

We consider several possibilities.  Suppose first that $|\gamma_1| = 1$, i.e.~$\gamma_1 = \{ i_1 \}$.  This then forces $\delta_2 = \{ i_3, i_5, \dotsc, i_{2m-1} \}$, $\gamma_2 = \{ i_2, i_4, \dotsc, i_{2m-2} \}$, and $\delta_1 = \{ i_{2m} \}$.   We can then factor $R$ as 
\begin{align*}
R &= \plucker{i_1,i_2,i_4,\ldots,i_{2m-2},\Delta}{i_3,i_5,\ldots,i_{2m-1},i_{2m},\Delta}{i_2,i_4,\ldots,i_{2m-2},i_{2m},\Delta}{i_1,i_3,i_5,\ldots,i_{2m-1},\Delta}  \\
&= \plucker{i_3,i_5,\ldots,i_{2m-1},i_{2m},\Delta}{i_1,i_3,i_4,i_6,\ldots,i_{2m-2},\Delta}{i_1,i_3,i_5,\ldots,i_{2m-1},\Delta}{i_3,i_4,i_6,\ldots,i_{2m-2},i_{2m},\Delta} \\
& \quad \cdot \plucker{i_3,i_4,i_6,\ldots,i_{2m-2},i_{2m},\Delta}{i_1,i_2,i_4,\ldots,i_{2m-2},\Delta}{i_1,i_3,i_4,i_6,\ldots,i_{2m-2},\Delta}{i_2,i_4,\ldots,i_{2m-2},i_{2m},\Delta},
\end{align*}
where the ellipses indicate the the sequence continues with the same parity subscripts.  Note that $\{i_1, i_5, i_7, \dotsc\}$ interlaces with $\{i_4, i_6, i_8, \dotsc\}$ and that $\{i_1, i_3\}$ interlaces with $\{i_2, i_4\}$, hence both ratios on the right hand side satisfy condition (M) by Remark~\ref{rem:interlacing_denominator_implies_M}.  Labeling the ratios on the right hand side as $R_1$ and $R_2$ respectively, we see that $\nu(R_1) = \nu(R) - 1 < \nu(R)$ and $\nu(R_2) = 2 < \nu(R)$.

If instead $|\gamma_1| > 1$, we may define define $k\geq 2$ to be the value so that $\{i_3,i_5,\dotsc,i_{2k-1}\} \subseteq \delta_2$ and $\{i_1,i_{2k+1}\}\subseteq \gamma_1$.
Similarly let $l \geq 1$ be the value so that $\{i_2,i_4,\ldots,i_{2l}\} \subseteq \gamma_2$ and $i_{2l+2} \in \delta_1$.  Observe that because both $\gamma_1$ interlaces with $\delta_1$ and $\gamma_2$ interlaces with $\delta_2$, we necessarily have $l=k-1$.

Let $\delta_{21} = \{i_3\}$, $\delta_{22} = \delta_2 \setminus \delta_{21}$, $\gamma_{21} = \{i_2\}$, $\gamma_{22} = \gamma_2 \setminus \gamma_{21}$.

We claim $\gamma_1 \cup \delta_{21}$ interlaces with $\delta_1 \cup \gamma_{21}$, as
\begin{align*}
	\gamma_1 \cup \delta_{21} &= \{ i_1, i_3, i_{2k+1}, \ast \} \text{ and} \\
	\delta_1 \cup \gamma_{21} &= \{ i_2, i_{2l+2}, \ast \},
\end{align*}
noting $2l+2 = 2k$ and $\gamma_1$ interlaces with $\delta_1$ by hypothesis.

Similarly, we claim $\gamma_1 \cup \delta_{22}$ interlaces with $\delta_1 \cup \gamma_{22}$, as
\begin{align*}
	\gamma_1 \cup \delta_{22} &= \{ i_1, i_5, i_7, \dotsc, i_{2k+1}, \ast \} \text{ and} \\
	\delta_1 \cup \gamma_{22} &= \{ i_4, i_6, \dotsc, i_{2l+2}, \ast \},
\end{align*}
noting $2l+2 = 2k$ and $\gamma_1$ interlaces with $\delta_1$ by hypothesis.

Therefore, we may write
\[
\plucker{\gamma_1,\gamma_2,\Delta}{\delta_1,\delta_2,\Delta}{\gamma_1,\delta_2,\Delta}{\delta_1,\gamma_2,\Delta}
=
\plucker{\gamma_1,\gamma_2,\Delta}{\delta_1,\delta_{21},\gamma_{22},\Delta}{\delta_1,\gamma_2,\Delta}{\gamma_1,\delta_{21},\gamma_{22},\Delta}
\cdot
\plucker{\gamma_1,\delta_{21},\gamma_{22},\Delta}{\delta_1,\delta_2,\Delta}{\delta_1,\delta_{21},\gamma_{22},\Delta}{\gamma_1,\delta_2,\Delta}
\]
where both ratios on the right hand side satisfy condition (M) by Remark~\ref{rem:interlacing_denominator_implies_M}.  Labeling the ratios on the right hand side as $R_1$ and $R_2$ respectively, we see that $\nu(R_1) = \nu(R) - |\gamma_{22}|$ and $\nu(R_2) = \nu(R) - |\delta_{21}|$.  Now $|\delta_{21}| = 1$, so $\nu(R_2) < \nu(R)$.

If $|\gamma_{22}| \geq 1$ we are finished.  If instead $|\gamma_{22}| = 0$, we deduce that
$\gamma_1 = \{i_1, i_5, i_7, \dotsc, i_{2m-1} \}$, $\delta_1 = \{ i_4, i_6, \dotsc, i_{2m} \}$, $\gamma_2 = \{i_2 \}$, and $\delta_2 = \{ i_3 \}$.  We can then factor $R$ as 
\begin{align*}
\plucker{i_1,i_2,i_5,i_7,\ldots,\Delta}{i_3,i_4,i_6,\ldots,\Delta}{i_2,i_4,i_6,\ldots,\Delta}{i_1,i_3,i_5,i_7,\ldots,\Delta} = &
\plucker{i_1,i_2,i_5,i_7,\ldots,\Delta}{i_3,i_5,i_6,i_8,\ldots,\Delta}{i_1,i_3,i_5,i_7,\ldots,\Delta}{i_2,i_5,i_6,i_8,\ldots,\Delta} \\
& \cdot \plucker{i_2,i_5,i_6,i_8,\ldots,\Delta}{i_3,i_4,i_6,\ldots,\Delta}{i_3,i_5,i_6,i_8,\ldots,\Delta}{i_2,i_4,i_6,\ldots,\Delta},
\end{align*}
where the ellipses indicate the the sequence continues with the same parity subscripts.  Note that $\{i_1, i_3, i_7, \dotsc\}$ interlaces with $\{i_2, i_6, i_8, \dotsc\}$ and that $\{i_3, i_5\}$ interlaces with $\{i_2, i_4\}$, hence both ratios on the right hand side satisfy condition (M) by Remark~\ref{rem:interlacing_denominator_implies_M}.  Labeling the ratios on the right hand side as $R_1$ and $R_2$ respectively, we see that $\nu(R_1) = \nu(R) - 1 < \nu(R)$ and $\nu(R_2) = 2 < \nu(R)$.

\end{pf}

\begin{thm}\label{thm:can_factor_if_nu(R)>3}
	Let $ R\in \mathcal B$ with $\nu(R) \geq 3$.  Then we may write $R=R_1R_2$ for some ratios $R_1,R_2\in \mathcal B$ with $\nu(R_i)<\nu(R)$ for $i=1,2$.
\end{thm}

\begin{pf}
If either $\gamma_1$ and $\delta_1$ do not interlace, or $\gamma_2$ and $\delta_2$ do not interlace, or both we may appeal to Lemma~\ref{lem:some_dont_interlace}.

If instead both $\gamma_1$ interlaces with $\delta_1$ and $\gamma_2$ interlaces with $\delta_2$, we may assume without loss of generality that $R$ is agreeably labeled (see Claim~\ref{claim:wlog_labeling}).  Observe that with this labeling $i_2 \in \delta_1\cup\gamma_2$.

If $i_2 \in \delta_1$, we may appeal to Lemma~\ref{lem:i_2-in-delta_1}.

If $i_2 \in \gamma_2$, we may appeal to Lemma~\ref{lem:i_2-in-gamma_2}.
\end{pf}

We now state our main result. 

\begin{thm}[Main Theorem]
Let $R$ be a ratio of the form $\alpha/\beta=\iplucker{\alpha_1}{\alpha_2}{\beta_1}{\beta_2}$ where $\alpha_1,\alpha_2,\beta_1,\beta_2$ are index sets in $\{1,\ldots,2n\}$.  The following are equivalent:
\begin{enumerate}
\item $R$ satisfies the ST0 condition and
\begin{equation}\max(| \alpha_1\cap
L|,| \alpha_2 \cap L|) \geq \max(| \beta_1\cap L|,| \beta_2\cap
L|)\label{secondcondition}\end{equation} for every interval
$L\subseteq \{1,\ldots,2n\},$ i.e. $\alpha/\beta$ satisfies condition (M).
\item $R$ can be written as a product of basic
ratios.
\item $R$ is bounded by 1.
\item $R$ is bounded.
\end{enumerate}\label{maintheorem}
\end{thm}

\begin{pf}
Note that $(2)\implies(3)\implies(4)\implies(1)$ are clear, so what
remains to show is that $(1)\implies(2)$.   

By Theorem~\ref{thm:can_factor_if_nu(R)>3}, we may write any ratio $R$ of the specified form satisfying condition (M) with $\nu(R) \geq 3$ as a product of ratios $R_1 R_2$ of the same form where each satisfies condition (M) and $\nu(R_i) < \nu(R)$ for $i=1,2$.
By Lemma~\ref{lem:nu(R)=2=>trivial-or-elementary} and the remarks directly preceding it, any ratio $R$ of the specified form satisfying condition (M) with $\nu(R) \leq 2$ is either trivial or can be written as a product of basic ratios.

A simple induction argument on the value of $\nu(R)$ completes the proof.
\end{pf}

\begin{rem}
Much of the proof in this section extends similar techniques used by Fallat et al.~to the case of non-principal minors (see \cite{Fallat:2003az}).
The equivalence of (1), (3) and (4) in Theorem \ref{maintheorem} is a result of Skandera  (see \cite{Skandera:2004qp}).
\end{rem}




\begin{rem}While we have shown that condition (M) implies boundedness for this specific class of ratios,  this condition is not sufficient in general.  For example, the ratio
\[\frac{\p{1,2,3,8}\p{2,3,4,5}\p{4,6,7,8}}{\p{1,4,6,8}\p{2,3,4,8}\p{2,3,5,7}}\] satisfies condition (M) but is unbounded over the class of totally positive matrices.  For example, when applied to the totally positive matrix
\[ \left(\begin{array}{llll}
1 & 3t^{-1} & 3t^{-2} & t^{-1} \\
2+t^{-1} & 1+ 6 t^{-1} + 3t^{-2} & 2t^{-1}+6t^{-2}+3t^{-3} & 1+2t^{-1} +t^{-2} \\
t+2 & t +4+6t^{-1} & 3+5t^{-1}+6t^{-2} & 2t+2+2t^{-1} \\
t & t+3 & t +2+3t^{-1} & t^2+t+2
\end{array}\right) \]
where $t$ is a positive indeterminate, the exhibited ratio increases without bound as $t\to\infty$.
\end{rem}

\begin{conj}
A ratio $\alpha / \beta$, where $\alpha$ and
$\beta$ are each sequences of an arbitrary number of index sets is bounded if and only if it can be written as a product of basic ratios.\label{conjecture3.5}
\end{conj}

This conjecture was briefly hinted at by Fallat et al.~with regards to a possible way to save a similar conjecture with respect to bounded ratios of principal minors (see \cite[\textsection 6]{Fallat:2003az}).

\section{Computational Methods and Computational Results}


Given this collection of basic ratios, a natural question to consider is: {\em  What is the set of ratios generated by products of positive powers of the basic ratios?}  Every ratio in this space is both bounded and expressible as a product of basic ratios.


We consider a typical element of this space to be a ratio of
products of index sets. Recall there are $N={2n \choose n}$ such
index sets. Then each ratio can be described by giving the power to
which each index set appears in the ratio; terms appearing in the
denominator have negative exponent. This allows us to identify each
ratio with a vector in $\R^N$ where each entry represents the power
to which that index set appears in the ratio.  The product of two
ratios then simply corresponds to the sum of their two associated
vectors in $\R^N$.

We write $v_1,\ldots,v_M$ as the vectors that correspond to each of
the $$M= n(2n-3) {2n-4 \choose n-2}$$ basic ratios.  Such a list of
generating vectors can be easily computed using {\tt Mathematica}.

The set of all ratios that are products of positive powers of
the basic ratios is a polyhedral cone \[P=\{ \sum_{i=1}^M \lambda_i v_i,
\,\, \lambda_i\in\R_{\geq 0}\}.\]
$P$ can also be
described as the intersection of finitely many linear half-spaces,
namely
\[ P = \{ x\in\R^N\, |\, A\cdot x \leq 0 \}\] for some matrix $A$ that can be computed. 
The software program {\tt cdd+} \cite{Fukuda:xb} is useful in the conversion between
these two representations of convex polyhedral cones.

We illustrate the utility of this by proving a non-trivial theorem.
This result was first obtained in \cite{Fallat:2000ps}.

\begin{prop}
Every bounded ratio of minors of a $3\times 3$ totally positive matrix can be
written as a product of positive powers of the basic ratios.
Furthermore, every such bounded ratio is bounded by
1.\label{basicratio}
\end{prop}

This can be verified computationally by computing the half planes of the cone generated by the basic ratios and then constructing a matrix in terms of a parameter $t$ that satisfies the inequality listed above.


As mentioned in \cite{Fallat:2003az}, one method of determining
whether or not a ratio is bounded is to work with a totally positive matrix
corresponding to the diagram in Figure \ref{fig:general}.  The entries in this matrix are then polynomials in the variables $l_i$, $d_j$, and $u_k$.  It is well known that this matrix will be totally positive if each variable is chosen to be positive, and all totally positive matrices may be arrived at by this construction for appropriate choices of the variables.

In this view, a ratio of minors  $R$ is a rational function $p/q$ in the same variables.  Some information about the ratio $p/q$ may be gleaned by examining the difference $q-p$ as a polynomial in the variables $l_i$, $d_j$, and $u_k$.  We denote this polynomial by $p_R$.  For example, if every coefficient in $p_R$ is positive (we call this `subtraction free') then $p_R$ will be positive for \emph{any} choice of positive variables $l_i$, $d_j$, and $u_k$.  This would imply that the ratio $p/q$ is necessarily bounded by 1 over the class of totally positive matrices.

This suggests the following conjecture, formulated in
\cite{Fallat:2003az,Skandera:2004qp}

\begin{conj}
A ratio of minors $R$ is bounded if and only if $p_R$ is subtraction free.\label{subfree}
\end{conj}
In other words, if $p_R$ is not subtraction
free, we conjecture that it is possible to find a family of 
totally positive matrices on which the ratio increases without bound.  (It has
\emph{always} been possible in every ratio that we have examined.)
This is significant, because of the existence of polynomials which
remain nonnegative but are not subtraction free. (e.g.
$x^2+y^2-2xy+1$).

\begin{rem}
Observe that Conjecture \ref{subfree} follows from Conjecture
\ref{conjecture3.5} and the short Pl\"ucker relation
(Equation~\ref{shortplucker}).  Indeed, the short Pl\"ucker relation guarantees that all basic ratios are subtraction free.  
This fact extends to arbitrary products of basic ratios, noting that if $R=A/B\cdot C/D$ where both $p_{A/B}$ and $p_{C/D}$ are
subtraction free, then $p_R = BD-AC=D(B-A)+A(D-C)$ is also
subtraction free.
\end{rem}

Using {\tt Mathematica}, we considered the set of ratios of the form
\begin{equation}\frac{\p{\alpha_1}\p{\alpha_2}\p{\alpha_3}}{\p{\beta_1}\p{\beta_2}\p{\beta_3}}\label{3over3}\end{equation}
over the class of $4\times 4$ totally positive matrices.  

Of the ratios satisfying the required ST0 and majorization conditions, we found that approximately 98\% could be written as a product of basic ratios. Those that could not be written as a product of basic ratios were found to be {\em not} subtraction free and actually unbounded over the class of totally positive matrices.
The results of these computer
experiments can be summarized in the following proposition.

\begin{prop}
For ratios of the form in Equation \ref{3over3}, the following are
equivalent when working over $4\times 4$ totally positive matrices:
\begin{enumerate}
\item The ratio is bounded;
\item The ratio
is bounded by 1;
\item The ratio can be factored into a product of
basic ratios; and
\item The ratio is subtraction free.
\end{enumerate}
\end{prop}

\section*{Acknowledgments}
The authors would like to thank their advisors, Misha Gekhtman and
Frank Connolly for the countless hours of help they provided.  The
authors also wish to thank the referee for numerous suggestions that
helped improve the quality of this paper.  This research was done
during a Notre Dame REU supported by NSF Grant DMS-0354132.

\bibliographystyle{elsart-num}
\bibliography{thesis}

\end{document}